\documentclass[10pt]{article}
\usepackage[utf8]{inputenc}
\usepackage[T1]{fontenc}
\usepackage[english]{babel}
\usepackage{amsmath, amsthm, amssymb,amsfonts}
\usepackage{array}
\usepackage{enumerate}
\usepackage[shortlabels]{enumitem}
\usepackage{mathrsfs}
\usepackage{mathtools}
\usepackage{theoremref}
\usepackage{doi}
\usepackage{hyperref}

\usepackage{stmaryrd} 
\usepackage[toc,page]{appendix}

\usepackage{dsfont}
\usepackage{enumitem}

\usepackage{thmtools}

\usepackage{algorithm}
\usepackage{algorithmic}

\usepackage{makeidx}
\makeindex

\usepackage{accents}

\usepackage{physics}
\usepackage{amsmath}
\usepackage{tikz}
\usepackage{mathdots}
\usepackage{yhmath}
\usepackage{cancel}
\usepackage{color}
\usepackage{siunitx}
\usepackage{array}
\usepackage{multirow}
\usepackage{amssymb}
\usepackage{gensymb}
\usepackage{tabularx}
\usepackage{extarrows}
\usepackage{booktabs}
\usetikzlibrary{fadings}
\usetikzlibrary{patterns}
\usetikzlibrary{shadows.blur}
\usetikzlibrary{shapes}

\usepackage{graphicx}

\usepackage{pgfplots}
\pgfplotsset{compat=newest}
\pgfplotsset{plot coordinates/math parser=false}
\usepackage{tikzscale}
\usepgfplotslibrary{external}
\usepackage{color}

\usepackage[official]{eurosym}
\usepackage{longtable}

\theoremstyle{definition}

\theoremstyle{remark}

\newcommand{\R}[0]{\mathbb{R}}

\newcommand{\cb}[0]{\mathrm{cb}}
\newcommand{\chp}[0]{\mathrm{chp}}
\newcommand{\hp}[0]{\mathrm{hp}}

\newcommand{\Esum}[0]{\mathrm{Esum}}
\newcommand{\Gsum}[0]{\mathrm{Gsum}}

\newcommand{\Emax}[0]{\mathrm{Emax}}
\newcommand{\Gmax}[0]{\mathrm{Gmax}}

\newcommand{\MEL}[0]{\mathrm{MAXELOAD}}
\newcommand{\SEL}[0]{\mathrm{SUMELOAD}}
\newcommand{\MGL}[0]{\mathrm{MAXGLOAD}}
\newcommand{\SGL}[0]{\mathrm{SUMGLOAD}}
\newcommand{\MHL}[0]{\mathrm{MAXHLOAD}}
\newcommand{\SHL}[0]{\mathrm{SUMHLOAD}}

\usepackage[backend=biber,style=alphabetic]{biblatex}
\title{An MINLP Model for designing decentralized energy supply networks}
\date{}
\author{Carl Eggen\thanks{Department of Mathematics and Statistics, University of Konstanz, 78457 Konstanz, Germany, {\texttt{\{carl.eggen, thanh-van.huynh, moritz.link, paul.stephan, stefan.volkwein\}@ uni-konstanz.de}}} \and Thanh-Van Huynh\footnotemark[1] \and Moritz Link\footnotemark[1] \and Paul Stephan\footnotemark[1] \and Stefan Volkwein\footnotemark[1]}

\begin{document}
\maketitle

\begin{abstract}
\noindent In this report, a detailed description of an MINLP model for decentralized energy supply network optimization is given. This model includes the possibility of extending gas transmission lines, local choice of heating technology, as well as local decisions for energy-efficient house renovation. Ultimately, the model is aimed at finding cost-efficient network plans while reducing carbon emissions to a specified amount.
\end{abstract}

\subsection*{Introduction}

In the following, we introduce and explain a mathematical model for designing decentralized energy supply networks. This model was developed by J.\ Lu together with D.\ König (Rechenzentrum für Versorgungsnetze Wehr GmbH, thereafter RZVN; see \url{https://www.rzvn.de}). Hence, if not marked otherwise, everything is based on the work from \cite[Chapter~3]{Lu2021} and \cite{Lu2018}.

The underlying energy supply concept of the model is based on small-scale decentralized technologies for heat and electricity generation, so-called micro-energy-conversion technologies (MECTs), which therefore provide the possibility for a (locally) more flexible mixture of energy resources. Practically speaking, energy consumers are able to decide individually which technology and therefore which energy carrier (electricity or (natural) gas) they use for covering their heating demand. This is where the decentralized aspect of the supply network comes into play, as the way in which heat is generated is decided at the local level. Additionally, the consumers are able to renovate their houses with the aim of energy efficiency, i.e.\ there is the possibility of a one-time investment in order to reduce future energy consumption (accompanied by lower energy consumption costs).

The mathematical model was designed to develop cost-efficient strategies for designing decentralized energy supply networks while simultaneously given targets for $\text{CO}_2$ emissions are achieved.

Since designing energy supply networks by hand is a very complex task -- especially with increasing size of the network -- it is necessary to use the power of mathematical modeling together with suitable optimization strategies. This report is aimed at providing and explaining such a mathematical model in all its details. Naturally, this model covers not all aspects of energy supply networks, but it is designed in such a way, that it can be used as a basis for decision-making in long-term energy policy.


Mathematically, the model is formulated as a \textit{Mixed-Integer Non-Li\-ne\-ar Pro\-gram} (MINLP). For an overview of this class of optimization problems, we refer the reader to \cite{Lee2012}. The nonlinearity comes from nonlinear equality constraints mimicking the physical laws of electricity and gas flow which are modeled explicitly. We will see the benefit of this explicit modeling later on. This explicit modeling guarantees that the proposed network plans are in line with, e.g., the capacity of the available gas pipes. More details regarding this issue can be found in the sections dealing with energy flows. The occurring discrete structure relies partially on binary decision variables representing possible extensions of transmission lines as well as allocation of MECTs. Furthermore, the modeling of energy-efficient house renovation using a piecewise linear approximation of a concave utility function also contributes to the discreteness of the problem.

This report is structured as follows: it starts with a description of the underlying network architecture (see p.~\pageref{ssec:network_archi}), after which the general idea of decentralized energy supply is introduced (see p.~\pageref{ssec:decentral_energy_supply}). Based on this, the modeling of energy-efficient house renovations is explained from page~\pageref{ssec:renovations} on. Afterward, the necessity, as well as the details regarding the nonlinear flow equations for electric current (see p.~\pageref{ssec:elec_supply}) and gas (see p.~\pageref{ssec:gas_supply}), are presented. For each energy carrier, an option of including the sizing of the respective transmission lines into the model is given (see p.~\pageref{sssec:cable_types} for electricity and p.~\pageref{sssec:pipe_types} for gas). Both are new optional extensions of the model which can be added individually. For the nonlinear gas flow, described by a simplification of the Darcy-Weisbach equation, we present a reformulation avoiding the $\mathrm{sign}$ function as proposed in \cite{BorrazSanchez2016} and \cite{Lu2021} (see p.~\pageref{ssec:sign_function}). From p.~\pageref{ssec:costs} on, all cost terms are explained and the objective function together with the constraints regarding the carbon emissions are stated. Finally, on p.~\pageref{ssec:variables}, the units of all appearing variables and parameters are summarized in a table. 

\subsection*{Network Architecture}
\label{ssec:network_archi}

The underlying structure of the model is a directed and connected graph $(V, E)$ consisting of the set of nodes $V$ and the set of arcs $E$. The set of nodes $V$ is given by $V = V_0 \cup V_1$, where $V_0 = \{0\}$ consists of the single source node and the set $V_1 = \{1,\ldots,n\}$ represents the remaining $n$ nodes. The energy for covering the demand of the whole network is injected at the source node and consumed at the arcs. The arcs $(i,j) \in E$ of the network can therefore be interpreted as streets of a city, i.e.\ an agglomeration of a certain number of residential houses, whereas the nodes can be seen as junctions. Thus, it is natural to assign the energy demand as well as the energy flows to the arcs. Consequently, the energy distribution takes place over the arcs, where we assume that the electricity transmission lines are already installed. There is an extension, where the elctricity transmission lines are not assumed to be a priorily installed. On the contrary, we do not assume the gas pipes have been constructed yet. Therefore, it is an important part of the optimization of the model to decide, where gas pipes need to be built in order to ensure the required distribution to the arcs. 

\subsection*{Decentralized Energy Supply}
\label{ssec:decentral_energy_supply}

\begin{figure}[h!]
	\centering
	\scalebox{.75}{\includegraphics[]{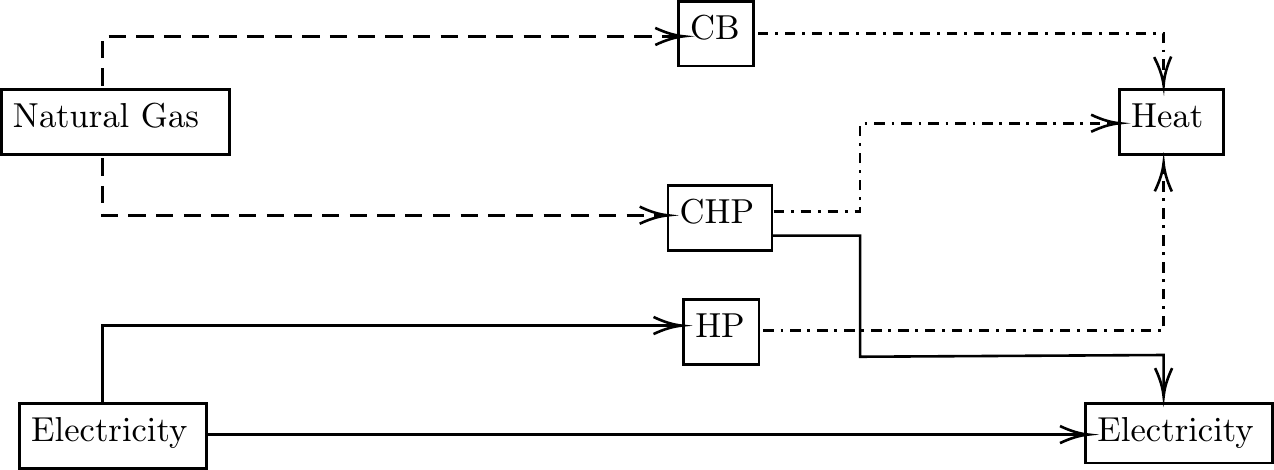}}
	\caption{Energy supply based on three MECTs (cf.\ \cite[Figure~3.1]{Lu2021})}
	\label{fig:MECTs}
\end{figure}

The main goal of the model is to meet the demand for electricity and heat at all arcs of the network. The yearly demand for electricity which is not used for heating, i.e.\ the amount of electricity to power, e.g., TV and refrigerator, is given by $\SEL_{i,j}$  (in kWh) for each arc $(i,j) \in E$, whereas the yearly demand for heat at arc $(i,j)$ is given by $\SHL_{i,j}$ (in kWh). Further, energy networks need to be designed in such a way, that the demand on peak times is met as well, since it has to be ensured that the capacity of the transmission lines is large enough to cover the maximal demands. Thus, the peak demand is given by $\MHL_{i,j}$ for gas (in kW) and $\MEL_{i,j}^t$ for electricity (in kW) for each arc $(i,j)$. In \Autoref{table:loads}, there is an example of the data for arc $(10,11)$ from one of our test networks. The data is coming from a simulation done by the project partners from RZVN. The simulation includes technical conditions of the MECTs together with energy-related conditions of the houses of the network.
\begin{table}[h]
	\centering
	\resizebox{\textwidth}{!}{
	\begin{tabular}{cccccc}
		\hline
		$(i,j)$ & $\MEL_{i,j}$ & $\SEL_{i,j}$ & $\MHL_{i,j}$ & $\SHL_{i,j}$  \\ 
		\hline
		$(10,11)$ &  $4.6$ & $14771$ & $47.6$ & $66116$ \\
		\hline
	\end{tabular}}
	\caption{An example of the simulated energy demand data of electricity and gas}
	\label{table:loads}
\end{table}
To meet the demand for heat, there are three different types of MECTs considered in the model, namely \textit{condensing boilers} (CB), \textit{combined heat and power units} (CHP) and \textit{heating pumps} (HP). As depicted in \Autoref{fig:MECTs}, \textit{condensing boilers} (CB) and \textit{combined heat and power units} (CHP) convert gas to heat, whereas \textit{heating pumps} (HP) convert electricity to heat. As by-product CHP also produces electricity and hence the gas and electricity supply cannot be considered independently.

In the model, we define $T \coloneqq \{\cb, \chp, \hp\} $  as the set of MECTs. At each arc $(i,j)\in E$ we introduce a binary variable
\begin{align*}
	x_{i,j}^\cb,\: x_{i,j}^\chp,\: x_{i,j}^\hp \in \{0,1\},
\end{align*}
for all technologies $t \in T$ indicating if MECT $t$ is installed on arc $(i,j)$, e.g., assume that HP is installed on arc $(i,j)\in E$ then $x_{i,j}^\hp = 1$, and $x_{i,j}^\hp = 0$ otherwise. Since MECTs are used for heating, we allow exactly one MECT only on those arcs having a positive heat demand. This leads to the constraint
\begin{align}
	\label{eq:atmost1MECT}
	x_{i,j}^\cb + x_{i,j}^\chp + x_{i,j}^\hp = 1 
\end{align}
for each arc $(i,j)\in E$ with $\SHL_{i,j} > 0$.

Since the different MECTs are powered by different energy carriers, the yearly demand for gas and electricity on each arc $(i,j) \in E$ depends on the installed technology and is given by the data $\SGL_{i,j}^t$ and $\SEL_{i,j}^t$ (both in kWh). Furthermore, energy networks have to be designed in such a way, that the demand on peak times can be met as well. For each arc $(i,j) \in E$ and each technology $t \in T$, this peak demand is represented by $\MGL_{i,j}^t$ for gas and $\MEL_{i,j}^t$ for electricity (both in kW).
Since CB is fueled by gas and only produces heat, the demand for electricity with CB installed equals the demand for electricity which is not used for heating, i.e.,
\begin{align*}
   \SEL^\cb_{i,j} = \SEL_{i,j} 
\end{align*}
for all arcs $(i,j) \in E$. Due to energy loss during conversion, the gas demand is higher than the heat demand in the case of CB and CHP, i.e.
\begin{align*}
   \SGL^t_{i,j} \geq \SHL_{i,j} 
\end{align*}
for each $(i,j) \in E$ and $t \in \{\cb,\chp \}$. The same properties hold for the peak demand. This can be also seen from the data in \Autoref{table:loads} and \Autoref{table:loadspertec}. Further, one can observe the decrease of electricity loads if CHP is installed compared to the case when CB is installed. On the other hand, the gas loads increase since additionally to heating, gas is also used to produce electricity. Next, there is no demand for gas at all if HP is installed, but a large increase in the demand for electricity, since HP only heats with electricity. With HP installed, the demand for energy in total is at its lowest compared to the other two technologies, because HP also uses thermal energy from the outside of a building. 

\begin{table}[h]
	\centering
	\resizebox{\textwidth}{!}{
	\begin{tabular}{cccccc}
		\hline
		$(i,j)$ & $t$ & $\MEL^t_{i,j}$ & $\SEL^t_{i,j}$ & $\MGL^t_{i,j}$ & $\SGL^t_{i,j}$  \\ 
		\hline
		$(10,11)$ & CB  &  $4.6$ & $14771$ & $52.9$ & $73478$ \\
		$(10,11)$ & CHP &  $4.3$ &  $5573$ & $55.8$ & $83588$ \\
		$(10,11)$ & HP  & $18.9$ & $36836$ &    $0$ &     $0$ \\
		\hline
	\end{tabular}}
	\caption{An example of the simulated energy demand data of electricity and gas per technology (cf.\ \cite[Table~3.1]{Lu2021}) \label{table:loadspertec}}
\end{table}

The amount of electricity $s_{i,j}^\Esum$ and gas $s_{i,j}^\Gsum$ which needs to be supplied on arc $(i,j)$ to cover the yearly demand for heat and electricity is given by the following simplified constraints 
\begin{align*}
	s^\Esum_{i,j} & = \sum_{t\in T} \SEL^t_{i,j} x^t_{i,j}, & s^\Gsum_{i,j} & = \sum_{t\in T} \SGL^t_{i,j} x^t_{i,j}.
\end{align*}
The same holds for the peak demands of electricity $s^\Emax_{i,j}$ and gas $s^\Gmax_{i,j}$:
\begin{align*}
	s^\Emax_{i,j} & = \sum_{t\in T} \MEL^t_{i,j} x^t_{i,j},&
	s^\Gmax_{i,j} & = \sum_{t\in T} \MGL^t_{i,j} x^t_{i,j}.
\end{align*}
Let us mention that these constraints become a bit more complex later due to the consideration of energy efficient house renovation.

Note that the variables $s_{i,j}^\Esum, s_{i,j}^\Emax$ might be negative in general if CHP is installed and produces more electricity than needed on this arc. However, in our data the parameters $\MEL^\chp_{i,j}$ und $\SEL^\chp_{i,j}$ appear non-negative and so all variables $s_{i,j}^\Esum$, $s_{i,j}^\Emax$, $s_{i,j}^\Gsum$ and $s_{i,j}^\Gmax$ are non-negative as well. 

The yearly and the peak demand are now defined on the arcs. This fits with their interpretation as streets lined with residential houses, where the energy is consumed. But for technical reasons, namely for balancing the power flows, it is necessary to artificially move the peak demand from the edges to the nodes. For each node $i \in V$ we therefore introduce two auxiliary variables $\bar{s}^\Emax_i, \bar{s}^\Gmax_i \in \R$ related to the peak demand of electricity and gas. The demand at the node $i$ is now defined as the half of all the demands at the arcs which are connected to node $i$:
\begin{align}
	\label{eq:artmove}
	\begin{aligned}
	    \bar{s}^\Emax_i&= \frac{1}{2} \left( \sum_{j\mid(j,i)\in E} s^\Emax_{j,i} + \sum_{k\mid(i,k)\in E} s^\Emax_{i,k} \right),\\
	    \bar{s}^\Gmax_i & = \frac{1}{2} \left( \sum_{j\mid(j,i)\in E} s^\Gmax_{j,i} + \sum_{k\mid(i,k)\in E} s^\Gmax_{i,k} \right).
	\end{aligned}
\end{align}

Lastly, for the single source node $0\in V$ we associate two variables $s_0^\Esum, s_0^\Gsum \in \R_0^+$ representing the amount of electricity and gas to be injected into the network in order to cover the yearly heat and electricity demand of the whole network. This is realized by the following constraints:
\begin{align}
	\label{eq:injconst}
	s_0^\Esum = \sum_{(i,j)\in E} s_{i,j}^\Esum\quad\text{and}\quad s_0^\Gsum = \sum_{(i,j)\in E} s_{i,j}^\Gsum.
\end{align}

\subsection*{Energy-Efficient House Renovation}
\label{ssec:renovations}

Energy-efficient house renovation is one of the key actions in the European Green Deal to target the climate neutrality of the EU in 2050, since about 40 \% of the energy consumption stems from buildings heating \footnote{European Commission (2019) What is the European Green Deal? Available at: \url{https://ec.europa.eu/commission/presscorner/detail/en/fs_19_6714} (Accessed: 08 Dec 2022).}. In order to reduce this consumption, we want to consider energy-efficient house renovation in our model as well, since an extensive reduction of heat demand can be achieved by executing certain energy-related house renovation, e.g., insulation or installing well-sealed windows.

\begin{figure}[h]
	\centering
    {\includegraphics[height=40mm, width=60mm]{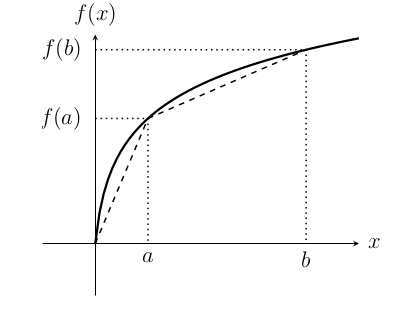}}
	\caption{Square root utility function of energy-efficient house renovation (black line) and its two-stage linear approximation (dashed line) (cf.\ \cite[Figure~3.2]{Lu2021})}
	\label{fig:utility}
\end{figure}

In the model this is realized by a two-stage technology-based approach, i.e.\ there are two renovation steps, which can only be realized one after another. For modeling the effect of the renovation investment mathematically, we use an utility function, which maps the investment costs to the reduction of the heat demand in percentage gained by this investment. We expect that one can significantly reduce the heat demand with comparable small investments in the beginning, while the gained reduction per invested euro decreases with progressing renovation. Thus we assume a concave utility function expressed by a square root function. The square root function is then approximated by a piecewise linear function with two straight-line segments, which yields the two-stage renovation approach as depicted in \Autoref{fig:utility}. As mentioned above, one can start with the second stage renovation only if the first stage one is already completed. Due to the square root shaped utility function the first stage renovation yields a larger amount of energy savings compared to the investment costs than the second stage renovations. 

In order to formulate the above model assumptions mathematically we associate $2\cdot \vert T\vert$ continuous and $\vert T\vert$ binary variables with each arc $(i,j)\in E$:
\begin{align*}
	x_{i,j}^{1,t} ,\: x_{i,j}^{2,t} \in [0,1] \textrm{ and } z^t_{i,j}\in \{0,1\} \textrm{ for all } t\in T.
\end{align*}
For $t\in T$, the variable $x_{i,j}^{1,t}$ stands for the progress of the first-stage renovation based on technology $t$\footnote{Clearly, one could reduce the total number of variables by using just one of each of the above variables for each arc instead of one for each technology and arc. This means, for example, that for each arc $(i,j)\in E$ we use the variable $x_{i,j}^1$ indicating to which extent the first-stage renovation is completed independently of the installed technology. However due to computational reasons it is more efficient the way it is done here.}, i.e.\ if, e.g., $x_{i,j}^{1,\mathrm{hp}} = 0.5$, it means that the first-stage renovation on arc $(i,j)$ with technology $\mathrm{HP}$  is $50\%$ done. Naturally, we can renovate on arc $(i,j)\in E$ based on $\mathrm{HP}$ only if $\mathrm{HP}$ is installed on arc $(i,j)$. This gives rise to the following inequality constraint:
\begin{align}
	\label{eq:installreno}
	x_{i,j}^t - x_{i,j}^{1,t} \geq 0 \quad \text{for all } (i,j)\in E,\: t \in T.
\end{align}
Furthermore, as explained before, starting with the second-stage renovation is possible only if the first-stage renovation is completed. This is expressed in the next inequality constraints for each $(i,j) \in E$:
\begin{align}
	\label{eq:firstsecren}
	x_{i,j}^{1,t} - x_{i,j}^{2,t} \geq 0, \quad
	x_{i,j}^{1,t} - z_{i,j}^t \geq 0, \quad
	x_{i,j}^{2,t} - z_{i,j}^t \leq 0
\end{align}
for all $t\in T$. In the above constraints the role of the binary variable $z_{i,j}^t$ becomes more clear. It serves as an indicator if the second-stage renovation is allowed to start\footnote{Note that one could also drop the binary indicator variable as in the optimization process finishing the first-stage renovation is always preferred to starting with the second-stage renovation due to the lower unit costs of the first-stage. However, it is not yet tested sufficiently to make a clear comment on the numerical comparison of these two options.}.

We proceed with the constraints regarding the goal of the renovations, namely the reduction of the energy demand. For a given arc $(i,j) \in E$ a completed renovation reduces the heat demand. As seen before, this affects also the gas and electricity demands of this arc. Therefore, we present updated versions of the constraints from the previous section. We start with the gas case, i.e.\ the constraints for $s_{i,j}^\Gsum$ and $s_{i,j}^\Gmax$:
\begin{align}
	\label{eq:gasrenored}
	\begin{split}
		s_{i,j}^\Gsum &= \sum_{t \in \{\cb,\chp\}} \SGL_{i,j}^t \left(x_{i,j}^t - \mu_1 x_{i,j}^{1,t} - \mu_2 x_{i,j}^{2,t} \right), \\
		s_{i,j}^\Gmax &= \sum_{t\in\{\cb,\chp\}} \MGL^t_{i,j} \left(x_{i,j}^t - \mu_1 x_{i,j}^{1,t} - \mu_2 x_{i,j}^{2,t} \right),
	\end{split}
\end{align}
where $\mu_1$ and $\mu_2$ represent the amount of energy savings in percentage of a fully completed renovation in the first stage and the second stage, respectively. Due to our two-stage modeling we have that 
\begin{align*}
	\mu_1 > \mu_2 > 0 \textrm{ and } \mu_1 + \mu_2 \leq 1.
\end{align*}

We continue with the electricity demand related to energy efficient house renovation. Here the situation is a bit more complicated than in the case of gas. This is due to the fact that renovation only affects the part of energy consumption which is related to heating. Depending on the technology $t\in T$ this part can be expressed by the following differences
\begin{align*}
    \SEL_{i,j}^t & - \SEL_{i,j}, \\
    \MEL_{i,j}^t &- \MEL_{i,j}.
\end{align*}
Consequently, the demand for electricity with the possibility of renovations is given by the constraints 
\begin{align}
    \label{eq:ereno}
	\begin{split}
		s_{i,j;t}^\Esum &= \SEL_{i,j} x_{i,j}^t \\ &+ \left(\SEL_{i,j}^t - \SEL_{i,j} \right) \left(x_{i,j}^t - \mu_1 x_{i,j}^{1,t} - \mu_2 x_{i,j}^{2,t}\right), \\
		s_{i,j;t}^\Emax &= \MEL_{i,j} x_{i,j}^t \\ &+ \left(\MEL_{i,j}^t - \MEL_{i,j} \right) \left(x_{i,j}^t - \mu_1 x_{i,j}^{1,t} - \mu_2 x_{i,j}^{2,t}\right),
	\end{split}
\end{align}
for each arc $(i,j) \in E$ and technology $t\in T$. Note that the first summand represents the demand for electricity independently of heating and the second summand corresponds to the demand for electricity only used for heating which is reduced by the renovation progress. If $\mathrm{HP}$ is installed on arc $(i,j)\in E$, we have an increased demand for electricity compared to the pure electricity demand, since electricity is also used for heating, i.e.\ $\SEL_{i,j}^\hp \geq \SEL_{i,j}$ and $\MEL^\hp_{i,j} \geq \MEL_{i,j}$. For CHP we recall that this technology also generates electricity as a by-product of the generation of heat from gas, i.e.\ if $\mathrm{CHP}$ is installed, the demand for electricity in total is less than the pure electricity demand, i.e.\ $\SEL_{i,j}^\chp \leq \SEL_{i,j}$ and $\MEL^\chp_{i,j} \leq \MEL_{i,j}$. Therefore, if renovation takes place, the amount of self-generated electricity decreases as less heat is needed and hence this yields an increase in the demand for electricity. In the case of CB, the constraints reduce to 
\begin{align*}
		s_{i,j;\cb}^\Esum = \SEL^\cb_{i,j} x_{i,j}^\cb\quad\text{and}\quad
		s_{i,j;\cb}^\Emax = \MEL^\cb_{i,j} x_{i,j}^\cb,
\end{align*}
since CB is not related to electricity, i.e.\ $\SEL_{i,j}^\cb = \SEL_{i,j}$ and $\MEL^\cb_{i,j} = \MEL_{i,j}$. 

Putting things together, for any arc $(i,j)\in E$ we obtain the reduced electricity demand
\begin{align}
	\label{eq:redelec}
	s^\Esum_{i,j} = \sum_{t\in T} s_{i,j;t}^\Esum \quad\text{and}\quad
	s^\Emax_{i,j} = \sum_{t\in T} s_{i,j;t}^\Emax.
\end{align}
%

\subsection*{Low-Voltage Electricity Supply}
\label{ssec:elec_supply}

In the following we describe how the electric power flow at peak time is modeled. We take this into account as any supply network needs to meet the demand at any time, but also has its own limitations, e.g., the public electricity supplier provides certain electric potential  $u_0$ at the source node $0 \in V_0$. As we will see, this initial potential together with the specific resistances of the underlying electric cables limit the electric power flow in the network. Thus, it may happen that for certain configurations we are not able to cover the demand for electricity on some parts of the network -- yielding in infeasibility of that specific configurations. This is not a problem in general, but it becomes a huge one if all possible configurations are infeasible. However, there are possibilities to overcome such inconveniences by, e.g., laying cables with less resistance or parallel cables in order to increase the possible electric power flow. Having this in mind, the importance of modeling (at least simplified) physical laws becomes apparent.

We proceed with the description of the modeling of the low-voltage electricity supply in its basic version. We associate two variables with the source node $0\in V_0$, namely the hourly peak amount of electricity required at the source node $s_0^\Emax$, which we already introduced, and the above mentioned new variable for the initial potential provided by the public electricity supplier 
\begin{align}
	\label{eq:sourcevolt}
	u_0 = u_{\max}.
\end{align}   
For each other node, i.e.\ for each $i\in V_1$, we denote by $u_i$ the electrical potential at that node. Also it must be guaranteed that at each node $i\in V_1$ the potential is neither higher than the initial potential nor lower than a minimal threshold potential:
\begin{align}
	\label{eq:voltbounds}
	0 \leq u_{\min} \leq u_i \leq u_{\max}.
\end{align}

As electric power flow takes place on the arcs $E$ of the network, for any arc $(i,j)\in E$ we introduce the variables $f_{\textrm{in},i,j}^\textrm{e}$ and $f_{\textrm{out},i,j}^\textrm{e}$ representing the power flow into $j$ and the power flow out of $i$, respectively. The electric power flows from $i$ to $j$ if the variables are positive and from $j$ to $i$ if they are negative. There are two different variables for the power flow on arc $(i,j)$ since there is a voltage drop on this arc and therefore $f_{\textrm{in},i,j}^\textrm{e} \neq f_{\textrm{out},i,j}^\textrm{e}$. The mentioned voltage drop on arc $(i,j)$ is represented by
\begin{align}
	\label{eq:voltdrop}
	\bar{u}_{i,j} = u_i - u_j.
\end{align}
\noindent

We are now ready to state the \emph{Ohmic law}, namely for any arc $(i,j)\in E$ we introduce the constraints:
\begin{align}
	\label{eq:ohmic}
	\begin{split}
		\textrm{R}^\textrm{e}_{i,j} f^\textrm{e}_{\textrm{in},i,j} &= a^\textrm{e} u_j \bar{u}_{i,j}, \\
		\textrm{R}^\textrm{e}_{i,j} f^\textrm{e}_{\textrm{out},i,j} &= a^\textrm{e} u_i \bar{u}_{i,j},
	\end{split}
\end{align}
where $a^\textrm{e}>0$ denotes the calorific multiplier of three-phase electric power flows and $\textrm{R}^\textrm{e}_{i,j}$ the resistance of the underlying electric cable. Furthermore, from \eqref{eq:ohmic} one immediately gets
\begin{align*}
	f^\textrm{e}_{\textrm{out},i,j} \geq f^\textrm{e}_{\textrm{in},i,j} 
\end{align*}
for all $(i,j) \in E$.
%

Lastly, for each node $i\in V$ we have to ensure that the electric power flow into node $i$ equals the sum of electric flow out of $i$ and the hourly peak demand for electricity $\bar{s}^\Emax_i$ at node $i$. This yields the following:
\begin{align}
	\label{eq:elecbalansink}
	\sum_{j\mid(j,i)\in E} f^\textrm{e}_{\textrm{in},j,i} - \sum_{k\mid(i,k)\in E} f^\textrm{e}_{\textrm{out},i,k} - \bar{s}_i^\Emax = 0,
\end{align}
for each node $i\in V_1$. The electric power flow into the source node $0\in V_0$ reduces to the hourly peak injection of electricity $s_0^\Emax$ (covering the hourly peak demands of the whole network) required at the source node. Therefore we get:
\begin{align}
	\label{eq:elecbalansource}
	s_0^\Emax - \sum_{k\mid(0,k)\in E} f^\textrm{e}_{\textrm{out},0,k} - \bar{s}^\Emax_0 = 0.
\end{align}
Putting \eqref{eq:elecbalansink} and \eqref{eq:elecbalansource} together we obtain the following:
\begin{align}
	\label{eq:injectmax}
	\begin{aligned}
	    s_0^\Emax - \sum_{i\in V} \bar{s}_i^\Emax \geq 0,
    \end{aligned}
\end{align}
This can be seen by using the above inequality for the inflow and outflow and reordering the summands in \eqref{eq:elecbalansink} and \eqref{eq:elecbalansource}:
\begin{align*}
    s_0^\Emax - \sum_{i\in V} \bar{s}_i^\Emax &=\sum_{k\mid(0,k)\in E} f^\textrm{e}_{\textrm{out},0,k}+\sum_{i\in V_1}\left(\sum_{k\mid(i,k)\in E} f^\textrm{e}_{\textrm{out},i,k}-\sum_{j\mid(j,i)\in E} f^\textrm{e}_{\textrm{in},j,i}\right)\\
    	&=\sum_{(i,j) \in E} \left( f^\textrm{e}_{\textrm{out},i,j} -  f^\textrm{e}_{\textrm{in},i,j}\right) \geq 0
\end{align*}
Let us mention that \eqref{eq:injectmax} is not part of the original model provided in \cite{Lu2021}. But it is added here since it enhances the numerical performance helping with the lower bounding.

In summary, the problem of electricity supply in our network can be formulated as follows: Assume that (given certain renovations) we have an hourly peak demand for electricity $\bar{s}^\Emax_i$ at each node $i\in V$. Are there feasible assignments for the node potentials $u_i$ and flows $f^\textrm{e}_{\textrm{in},i,j}$ and $f^\textrm{e}_{\textrm{out},i,j}$ such that the balance equations at each node together with the Ohmic law at each arc are satisfied?

\subsubsection*{Option: different cable types}	
\label{sssec:cable_types}

As mentioned before, it might be necessary to consider different cable types in order to ensure feasibility of the model. In order to do so, we use a similar approach as presented in \cite[Section~2.3]{Belotti2013} as a possible extension of the original model from \cite{Lu2021}. Therefore, assume that $r$ cable types with different resistances are available. For any $k \in \{1,\ldots,r\}$ let $y^{\text{e}}_{i,j;k} \in \{0,1\}$ be the variable that indicates whether a cable of type $k$ is laid at edge $(i,j)$. We require that
\begin{align*}
	\sum_{k=1}^r y^{\text{e}}_{i,j;k} = 1	
\end{align*}
holds for all $(i,j) \in E$, i.e.\ we want to ensure that exactly one cable is laid at each arc.

Let $R_k^\text{e}$ be the resistance of a cable of type $k \in \{1,\ldots,r\}$ per meter of cable and $l_{i,j} \geq 0$ be the length of arc $(i,j)$ (in m). We define the resistance constant for each arc $(i,j) \in E$ by
\begin{align*}
	R_{i,j}^\text{e} := \sum \limits_{k=1}^r y^{\text{e}}_{i,j;k} R_k^\text{e} l_{i,j}.	
\end{align*}

Further, we denote by $\zeta^\text{e}_k$ the cost per meter for lying a cable of type $k$. Thus, the installation costs $\zeta_{i,j}^\text{e}$ of the chosen cable type on arc $(i,j)$ are given by 
\begin{align*}
	\zeta_{i,j}^\text{e} = \sum_{k=1}^r \zeta_k^\text{e} \: y^{\text{e}}_{i,j;k} \: l_{i,j}.
\end{align*}
In contrast to the original model, we do not assume that the electricity transmission lines are a-priorily available and therefore we need to invest in their installation. The costs for providing a suitable electricity grid for the network are then given by
\begin{align*}
	C_\text{grid}^\text{e} = \sum_{(i,j) \in E} \zeta^\text{e}_{i,j}.
\end{align*}
%

\subsection*{Low-Pressure Gas Supply}
\label{ssec:gas_supply}

The modeling of the gas supply is similar to the above one for the electricity supply. But there are some significant differences. Firstly, we assume that there does not exist a gas network yet, i.e.\ the optimization needs to decide on which arc gas pipes have to be laid. This decision interferes with the decision of choosing the technology for an arc since CB and CHP are fueled by gas. Secondly, the physics of gas flows differ from the ones of electric current. There is pressure instead of electrical potential and the Ohmic law is replaced by a simplification of the \emph{Darcy-Weisbach equation} for low-pressure gas flows.

For the source node, we have the variable $p_0$ representing the gas pressure provided by the public gas supplier, i.e.
\begin{align}
	\label{eq:gpsource}
	p_0 = p_{\max}
\end{align}
as well as the hourly peak injection at the source node $s_0^\Gmax$.

The potential at each node turns into the gas pressure $p_i$ at each node $i\in V$, and again, it has to be ensured that the pressure at each node is neither higher than the provided maximal pressure nor lower than a minimal threshold, i.e.
\begin{align}
	\label{eq:gpbounds}
	0 \leq p_{\min} \leq p_i \leq p_{\max}.
\end{align}

We recall that we do not assume the gas grid to be built a-priorily, i.e.\ the optimization of the model chooses the arcs where gas pipes need to be built. Therefore we associate two variables with each arc $(i,j)\in E$, namely the binary decision variable $y^g_{i,j}$ for building a gas pipe on arc $(i,j)$ as well as $f^\textrm{g}_{i,j}$ representing the gas flow\footnote{Again, $f^\textrm{g}_{i,j}$ is positive if gas flows from $i$ to $j$, and negative otherwise.} on arc $(i,j)$. Naturally, if there is positive gas demand on an arc $(i,j)\in E$, i.e.\ if $\mathrm{CB}$ or $\mathrm{CHP}$ are installed at that arc, it has to be linked to the gas network. This yields the following two constraints:
\begin{align}
	\label{eq:pipebuild}		x_{i,j}^\cb - y^g_{i,j} \leq 0\quad\text{and}\quad x_{i,j}^\chp - y^g_{i,j} \leq 0.
\end{align}
Furthermore, the pressure loss on an arc $(i,j)\in E$ is represented by
\begin{align}
	\label{eq:gploss}
	\bar{p}_{i,j} = p_i - p_j.
\end{align}
If on a given arc $(i,j)\in E$ a gas pipe has been built, i.e.\ if $y^g_{i,j} = 1$ holds, then the announced simplification of the Darcy-Weisbach equation comes into play:
\begin{align}
	\label{eq:DWeq}
	\big(a^\textrm{g}\big)^2 \bar{p}_{i,j} y^g_{i,j} = \textrm{R}^\textrm{g}_{i,j} \textrm{ sign} \big(\bar{q}_{i,j}\big) \big(f^\textrm{g}_{i,j}\big)^2,
\end{align}
where $a^\textrm{g}>0$ is the calorific multiplier of gas flow and $\textrm{R}^\textrm{g}_{i,j}>0$ is the resistance of the underlying gas pipeline. The variable $\bar{q}_{i,j}$ represents the amount of gas flow on arc $(i,j)\in E$ which works as an indicator of the direction of gas flow via the $\textrm{sign}$ function in \eqref{eq:DWeq} as we can see by the following:
\begin{align*}
    \bar{p}_{i,j}\left\{
    \begin{aligned}
    &\ge0&\text{if }\bar{q}_{i,j}\ge 0,\\
    &\le0&\text{if }\bar{q}_{i,j}\le 0.
    \end{aligned}
    \right.
\end{align*}
We also add\footnote{Note that this constraint has not been  part of the original model in \cite{Lu2021}.} the following constraint for any arc $(i,j)\in E$:
\begin{align}
	\label{eq:builgaspiperel}
	f^\textrm{g}_{i,j} = y^g_{i,j} f^\textrm{g}_{i,j},
\end{align}
which is redundant by \eqref{eq:DWeq}. But if we use a linearization of \eqref{eq:DWeq}, the above constraint may come into play and may yield better lower bounding. \hfill\\
\noindent
Again, similar as in the electricity setting, for each node $i\in V$ we have the balance equations of the power flow. For the nodes $i\in V_1$ we have:
\begin{align}
	\label{eq:gasbalansink}
	\sum_{j\mid(j,i)\in E} f^\textrm{g}_{j,i} - \sum_{k\mid(i,k)\in E} f^\textrm{g}_{i,k} - \bar{s}^\Gmax_i = 0.
\end{align}
And for the source node $0\in V_0$ the inflow reduces to the hourly peak injection of gas $s_0^\Gmax$ required at the source node. Consequently, we obtain the following constraint:
\begin{align}
	\label{eq:gasbalansource}
	s_0^\Gmax - \sum_{k\mid(0,k)\in E} f^\textrm{g}_{0,k} - \bar{s}^\Gmax_0 = 0.
\end{align}

Let us mention that we could add an inequality constraint similar to \eqref{eq:injectmax}. However, it turns out that this additional constraint is numerically not necessary.

\subsubsection*{Option: different pipes}
\label{sssec:pipe_types}

As mentioned above for the electricity supply, it might be necessary to consider different pipes in order to ensure the feasibility of the model. These pipes differ in diameter. A larger diameter means that more gas can flow through the pipe while having the same pressure. Assume that $r$ different pipe types with corresponding diameters $d_k$, $k\in\{1,\ldots,r\}$, are available. For any $k \in \{1,\ldots,r\}$ let $y^{\text{g}}_{i,j;k} \in \{0,1\}$ be a variable indicating whether a pipe of type $k$ is built at arc $(i,j)$. Again we ensure that only one pipe is laid at arc $(i,j)$ -- if gas pipes are required at all at that arc -- with the following constraint
\begin{align*}
	\sum_{k=1}^r y^{\text{g}}_{i,j;k} = y^g_{i,j}.	
\end{align*}
%
%
We denote by $\zeta^\text{g}_k$ the cost per meter for lying a pipe of type $k$. Thus, the cost $\zeta_{i,j}^\text{g}$ of laying the chosen pipe at arc $(i,j)\in E$ is given by 
\begin{align*}
	\zeta_{i,j}^\text{g} = \sum_{k=1}^r \zeta_k^\text{g} \: y^{\text{g}}_{i,j;k} \: l_{i,j}.
\end{align*}

Finally, we compute the resistance constant $R_{i,j;k}^\text{g}$ for pipe type $k$ on arc $(i,j)$ by the following formula as done in \cite{Zelmer2010}:
\begin{align*}
	R^{\text{g}}_{i,j;k} = \lambda\,\frac{8 \rho l_{i,j}}{\pi^2 d_k^5},
\end{align*}
where $\lambda = 0.3164/(v\rho d_k)^{0.25}$ is the \emph{Darcy friction coefficient}\footnote{We neglect the dynamic viscosity and therefore approximate the kinetic viscosity $\nu$ by $\nu \approx 1/\rho$.}, $\rho = 0.7\,\text{kg}/\text{m}^3$ the density and $v = 6\,\text{m}/\text{s}$ the velocity of the gas. Consequently, we obtain
\begin{align*}
	R_{i,j}^\text{g} := \sum_{k=1}^r y^{\text{g}}_{i,j;k} R^\text{g}_{i,j;k}.	
\end{align*}

\subsection*{Reformulation avoiding $\textrm{sign}$ function}
\label{ssec:sign_function}

In \eqref{eq:DWeq} the $\textrm{sign}$ function occurs. As most solvers cannot deal with it, the corresponding constraints are reformulated such that we can avoid the $\textrm{sign}$ function. In order to do so, an approach presented in \cite{BorrazSanchez2016} is applied.

For any arc $(i,j)\in E$ there are two binary variables $y_{i,j}^{g+}$ and $y_{i,j}^{g-}$ representing the gas flow direction:
\begin{align}
	\label{eq:refgflow}
	\begin{split}
		\big(1 - y_{i,j}^{g+}\big) \bar{q}_{\min} \leq \bar{q}_{i,j} &\leq \big(1 - y_{i,j}^{g-}\big) \bar{q}_{\max}, \\
		\big(1 - y_{i,j}^{g+}\big) \bar{p}_{\min} \leq \bar{p}_{i,j} &\leq \big(1 - y_{i,j}^{g-}\big) \bar{p}_{\max}, \\
		y_{i,j}^{g+} + y_{i,j}^{g-} &= 1,
	\end{split}
\end{align}
where 
\begin{align*}
	\begin{split}
		\bar{q}_{\min} &:= -\bar{q}_{\max} \leq \bar{q}_{i,j} \leq \bar{q}_{\max}, \\
		\bar{p}_{\min} &:= -\bar{p}_{\max} \leq \bar{p}_{i,j} \leq \bar{p}_{\max} = p_{\max} - p_{\min},	
	\end{split}
\end{align*}
and $\bar{q}_{\max}>0$ the maximal positive amount of gas flow allowed in the network. By the above constraints it is ensured that $\bar{q}_{i,j} \geq 0$ if and only if $y_{i,j}^{g+} = 1$ and therefore $\bar{q}_{i,j}$ has the same sign as $\bar{p}_{i,j}$. Thus, we can rewrite \eqref{eq:DWeq} without the $\textrm{sign}$ function as
\begin{align}
	\label{eq:newDWeq}
 	\big(a^\textrm{g}\big)^2 \big(y_{i,j}^{g+} - y_{i,j}^{g-}\big) \bar{p}_{i,j} y^g_{i,j} = \textrm{R}^\textrm{g}_{i,j} \big(f^\textrm{g}_{i,j}\big)^2.
\end{align}

We now want to linearize the nonlinear term on the left-hand side of \eqref{eq:newDWeq}. We start by introducing placeholder variables $\bar{p}_{i,j}^+$:
\begin{align}
	\label{eq:pbarplace}
	\bar{p}_{i,j}^+ = \big(y_{i,j}^{g+} - y_{i,j}^{g-}\big) \bar{p}_{i,j}
\end{align}
and then write \eqref{eq:newDWeq} as
\begin{align}
	\label{eq:DWeqplace}
	\big(a^\textrm{g}\big)^2 \bar{p}_{i,j}^+ y^g_{i,j} = \textrm{R}^\textrm{g}_{i,j} \big(f_{i,j}^\textrm{g}\big)^2
\end{align}
with $0 \leq \bar{p}_{i,j}^+ \leq \bar{p}_{\max}$. The remaining nonlinear term on the left hand side of \eqref{eq:DWeqplace} now only consists of the product of a continuous and a binary variable. Using McCormick relaxations \cite{McCormick1976} as done in \cite[Remark~2.3.2]{Lu2021} we can linearize \eqref{eq:pbarplace} as:
\begin{align}
	\label{eq:pbarref}
	\begin{split}
		\bar{p}^+_{i,j} &\geq \big(y_{i,j}^{g+} - y_{i,j}^{g-}\big) \bar{p}_{\max} + \bar{p}_{i,j} - \bar{p}_{\max}, \\
		\bar{p}^+_{i,j} &\geq \big(y_{i,j}^{g+} - y_{i,j}^{g-}\big) \bar{p}_{\min} - \bar{p}_{i,j} + \bar{p}_{\min}, \\
		\bar{p}_{i,j}^+ &\leq \big(y_{i,j}^{g+} - y_{i,j}^{g-}\big) \bar{p}_{\min} + \bar{p}_{i,j} - \bar{p}_{\min}, \\
		\bar{p}_{i,j}^+ &\leq \big(y_{i,j}^{g+} - y_{i,j}^{g-}\big) \bar{p}_{\max} - \bar{p}_{i,j} + \bar{p}_{\max}, \\
	\end{split}
\end{align}
and \eqref{eq:DWeqplace} as:
\begin{align}
	\label{eq:DWfinal}
	\begin{split}
		\textrm{R}^\textrm{g}_{i,j} \big(f^\textrm{g}_{i,j}\big)^2 &\geq \big(a^\textrm{g}\big)^2 \big(\bar{p}^+_{i,j} + \bar{p}_{\max} y^g_{i,j} - \bar{p}_{\max}\big), \\
		\textrm{R}^\textrm{g}_{i,j} \big(f^\textrm{g}_{i,j}\big)^2 &\leq \big(a^\textrm{g}\big)^2 \bar{p}^+_{i,j}, \\
		\textrm{R}^\textrm{g}_{i,j} \big(f^\textrm{g}_{i,j}\big)^2 &\leq \big(a^\textrm{g}\big)^2 \bar{p}_{\max} y^g_{i,j}.  
	\end{split}
\end{align}
Now we can replace the original simplification of the Darcy-Weisbach equation \eqref{eq:DWeq} by \eqref{eq:refgflow}, \eqref{eq:pbarref} and \eqref{eq:DWfinal}, where te remaining nonlinearity is only the quadratic term on the left hand side of \eqref{eq:DWfinal}.

\subsection*{Cost Minimization}
\label{ssec:costs}

The model is designed to provide possible network plans such that the demands at every arc are covered. The optimization now aims at finding the network plans that minimize the following costs:
\begin{itemize}
	\item operational costs: energy consumption costs to cover the yearly energy demand, energy allocation costs to cover peak demand,
	\item investment and maintenance costs: annualized costs related to investments in gas pipelines, MECTs and house renovations over a time horizon of 20 years.
\end{itemize}

We start with the operational costs. For energy consumption costs we consider the purchasing cost for electricity and gas:
\begin{align}
	\label{eq:enconpur}
	C_\textrm{energy} = \alpha_\textrm{p}^\textrm{e} s_0^\Esum + \alpha_\textrm{p}^\textrm{g} s_0^\Gsum,
\end{align}
where $\alpha_\textrm{p}^\textrm{e}>0$ and $\alpha_\textrm{p}^\textrm{g}>0$ are the before-tax unit prices in $\text{\euro{}}/ \textrm{kWh}$ of electricity and gas, respectively. Naturally, we have to add the tax costs for energy:
\begin{align}
	\label{eq:encontax}
	C_\textrm{tax} = \beta^\textrm{e} \left(s_0^\Esum + t_\textrm{adv} s_\chp^\Esum\right) + \beta^\textrm{g} s_0^\Gsum,
\end{align}
where $\beta^\textrm{e}>0$ and $\beta^\textrm{g}>0$ are the tax rates for electricity and gas per kWh. The number $s_\chp^\Esum$ is the total amount of electricity produced locally by the $\textrm{CHPs}$ and therefore has to be taxed as well. As self-production should be promoted, a tax advantage factor $t_\textrm{adv} \in (0,1)$ is applied.

The allocation costs depend on the peak injections and we get:
\begin{align}
	\label{eq:allocos}
	C_\textrm{allocation} = \alpha_\textrm{a}^\textrm{e}s_0^\Emax + \alpha^\textrm{g}_\textrm{a} s_0^\Gmax, 
\end{align}
where $\alpha_\textrm{a}^\textrm{e}>0$ and $\alpha_\textrm{a}^\textrm{g}>0$ are unit costs for electricity and gas per year in $\text{\euro{}}/ \textrm{kWa}$. This completes the operational costs.

We proceed with the investment and maintenance part of the total costs. We have seen that there is the possibility of building a gas network, i.e.\ the cost $\zeta^\textrm{g}_{i,j}$ for laying a gas pipeline on arc $(i,j)\in E$ in \euro{}/a leads to the annualized gas grid costs:
\begin{align}
	\label{eq:gridcosts}
	C_\textrm{grid} = \sum_{(i,j)\in E} y^g_{i,j} \zeta^\textrm{g}_{i,j}.
\end{align}
The annualized installation costs of the MECTs are obtained by
\begin{align}
	\label{eq:techcosts}
	C_\textrm{tech} = \sum_{t\in T} \sum_{(i,j)\in E} \gamma^t x^t_{i,j},
\end{align}
with $\gamma^t>0$ being the investment and maintenance cost for technology $t$ in \euro{}a.

Finally, we need to consider the annualized house renovation costs. The unit costs of first- and second-stage renovation $\nu_1$ and $\nu_2>0$, respectively, are provided in $\text{\euro{}}/\textrm{kWh}$ measured against the heat demand $\textrm{SUMHLOAD}_{i,j}$ on arc $(i,j)\in E$. Thus, the total renovation costs are given by
\begin{align}
	\label{eq:renocosts}
	C_\textrm{renov} = \sum_{t\in T} \sum_{(i,j) \in E} \textrm{SUMHLOAD}_{i,j} \left(\nu_1 x_{i,j}^{1,t} + \nu_2 x_{i,j}^{2,t} \right).
\end{align}

Another goal of the modeling is the limitation of $\text{CO}_2$ emissions in the network up to a given target, i.e.
\begin{align}
	\label{eq:carbontar}
	E_\textrm{carbon} - E_\textrm{target} \leq 0,
\end{align}
where the $\text{CO}_2$ emissions of the network are computed as follows:
\begin{align}
	\label{eq:carbon}
	E_\textrm{carbon} = \kappa^\textrm{e} s_0^\Esum + \kappa^\textrm{g} s_0^\Gsum
\end{align}
with $\kappa^\textrm{e}>0$ and $\kappa^\textrm{g}>0$ denoting the amount of $\text{CO}_2$ emissions in $\textrm{kg}/\textrm{kWh}$ of electricity and gas, respectively. \\
We are now ready to state the complete model:
\begin{align}
	\label{eq:model}
	\tag{Network}
	\left\{
	\begin{aligned}
	    &\min C= C_\textrm{energy} + C_\textrm{tax} + C_\textrm{allocoation} + C_\textrm{grid} + C_\textrm{tech} + C_\textrm{renov}\\
		&\hspace{0.5mm}\text{subject to (s.t.) the constraints \eqref{eq:atmost1MECT}-\eqref{eq:carbon}.}
	\end{aligned}
	\right.
\end{align}
\newpage
\subsection*{Variables and parameters}
\label{ssec:variables}

		\begin{longtable}{l l p{7cm}}
			\hline \\
			Sets					& $V$				& Nodes				\\
			\:						& $V_0$				& Source node		\\
			\:						& $V_1$				& Nodes without source node \\
			\:						& $E$				& Arcs				\\
			\:						& $T$				& technologies				\\

			Cost variables			& $C_\text{energy}$	& Overall energy cost of gas and electricity [\euro/a] \\
			\:						& $C_\text{tax}$	& Overall tax cost for gas and electricity [\euro/a] \\
			\:						& $C_\text{allocation}$	& Overall allocation cost of gas and electricity [\euro/a] \\
			\:						& $C_\text{grid}^\text{g}$	& Overall cost for installing the gas grid [\euro/a] \\
			\:						& $C_\text{tech}$	& Overall investment and maintenance cost for installed technologies [\euro/a] \\
			\:						& $C_\text{renov}$	& Overall renovation cost [\euro/a] \\
			\:						& $E_\text{carbon}$	& Overall CO2 emission of the network [kg/a] \\ 
			\:						& $C$				& Overall cost [\euro/a] \\
			\:						& $C_\text{grid}^\text{e}$	& Overall cost of installing electricity grid [\euro/a] \\

			Network variables		& $u_i$				& Electric potential at node $i$ [V] \\
			\:						& $\bar{u}_{i,j}$	& Electric voltage drop on arc $(i,j)$ [V] \\
			\:						& $f^\textrm{e}_{\textrm{in},i,j}$	& Electric power flow into node $j$ coming from arc $(i,j)$ [kW] \\
			\:						& $f^\textrm{e}_{\textrm{out},i,j}$	& Electric power flow out of node $i$ on arc $(i,j)$ [kW] \\
			\:						& $p_i$				& Gas pressure at node $i$ [mbar] \\
			\:						& $f_{i,j}^\textrm{g}$				& Gas power flow on arc $(i,j)$ [kW] \\ 
			\:						& $\bar{p}_{i,j}$	& Gas pressure loss on arc $(i,j)$ [mbar] \\
			\:						& $\bar{q}_{i,j}$	& Directed amount of gas flow on arc $(i,j)$ [m$^3$/h] \\
			\:						& $\bar{p}^+_{i,j}$	& Nonnegative placeholder variable for pressure loss on arc $(i,j)$ [mbar] \\
			\:						& $l_{i,j}$			& Length of arc $(i,j)$ [m] \\

			Binary variables		& $x^t_{i,j}$		& Binary decision for installing technology $t$ on arc $(i,j)$ \\
			\:						& $z^t_{i,j}$		& Binary decision for completing first stage reno\-va\-tion on arc $(i,j)$ w.r.t.\ technology $t$ \\
			\:						& $y^g_{i,j}$			& Binary decision for installing a gas pipe on arc $(i,j)$ \\	
			\:						& $y^+_{i,j}$		& Binary decision for gas flow from node $i$ to node $j$ \\
			\:						& $y^-_{i,j}$		& Binary decision for gas flow from node $j$ to $i$ \\
			\:						& $y_{i,j;k}$		& Binary decision for installing cable type $k$ on arc $(i,j)$ \\



			Hub system variables	& $s^\Emax_{i,j}$	& Hourly peak amount of electricity	supplied to arc $(i,j)$ [kW]	\\
			\:						& $s^\Esum_{i,j}$	& Yearly amount of electricity supplied to arc $(i,j)$ [kWh/a] \\
			\:						& $s^\Gmax_{i,j}$	& Hourly peak amount of gas supplied to arc $(i,j)$ [kW] \\
			\:						& $s^\Gsum_{i,j}$	& Yearly amount of gas supplied to arc $(i,j)$ [kWh/a] \\
			\:						& $\bar{s}_i^\Emax$	& Hourly peak amount of electricity supplied to node $i$ [kW] \\
			\:						& $\bar{s}_i^\Gmax$ & Hourly peak amount of gas supplied to node $i$ [kW] \\
			\:						& $s^\Emax_0$		& Hourly peak amount of electricity injected at the source node [kW] \\
			\:						& $s^\Esum_0$		& Yearly amount of electricity injected to the source node [kWh/a] \\
			\:						& $s^\Gmax_0$		& Hourly peak amount of gas injected to the source node [kW] \\
			\:						& $s^\Gsum_0$		& Yearly amount of gas injected to the source node [kWh/a] \\
			\:						& $x_{i,j}^{1,t}$	& Progress of first stage renovation on arc $(i,j)$ w.r.t.\ technology $t$ \\
			\:						& $x_{i,j}^{2,t}$	& Progress of second stage renovation on arc $(i,j)$ w.r.t.\ technology $t$ \\
			\:						& $s_{i,j;t}^\Emax$	& Hourly peak amount of electricity supplied to arc $(i,j)$ w.r.t.\ technology $t$ [kW] \\
			\:						& $s_{i,j;t}^\Esum$	& Yearly amount of electricity supplied to arc $(i,j)$ w.r.t.\ technology $t$ [kWh/a] \\
			\:						& $s_\chp^\text{Esum}$	& Yearly amount of electricity locally produced by CHPs [kWh/a] \\

			Hub system parameters	& $\MEL^t_{i,j}$	& Hourly peak demand for electricity on arc $(i,j)$ if technology $t$ is installed [kW]\\
			\:						& $\SEL^t_{i,j}$	& Yearly demand for electricity on arc $(i,j)$ if technology $t$ is installed [kWh/a] \\
			\:						& $\MGL^t_{i,j}$	& Hourly peak demand for gas on arc $(i,j)$ if technology $t$ is installed [kW] \\
			\:						& $\SGL^t_{i,j}$	& Yearly demand for gas on arc $(i,j)$ if technology $t$ is installed [kWh/a] \\
			\:						& $\mu_1$			& Percentage amount of energy savings of completed first stage renovation \\
			\:						& $\mu_2$			& Percentage amount of energy savings of completed second stage renovation \\
			\:						& $\text{SUMHLOAD}_{i,j}$	& Yearly heat demand on arc $(i,j)$ [kWh/a] \\



			Network parameters		& $u_{\max}$			& Electric potential provided by the public electricity supplier [V] \\
			\:						& $u_{\min}$			& Minimal threshold for electric potential [V] \\
			\:						& $a^\textrm{e}$		& Calorific multiplier of three--phase electric power flow \\
			\:						& $R_{i,j}^\textrm{e}$	& Electric resistance on arc $(i,j)$ [$\Omega$] \\
			\:						& $p_{\max}$			& Gas pressure provided by the public gas supplier [mbar] \\
			\:						& $p_{\min}$			& Minimal threshold for gas pressure [mbar] \\
			\:						& $a^\textrm{g}$		& Calorific multiplier of gas flow [J/m$^3$] \\
			\:						& $R_{i,j}^\textrm{g}$	& Gas resistance on arc $(i,j)$ [mbar/(m$^3$/h)$^2$] \\
			\:						& $\bar{q}_{\max}$		& Maximal positive amount of gas flow allowed in the network [m$^3$/h] \\
			\:						& $\bar{q}_{\min}$		& Maximal negative amount of gas flow allowed in the network [m$^3$/h] \\
			\:						& $\bar{p}_{\max}$		& Maximal positive pressure loss allowed in the network [mbar] \\
			\:						& $\bar{p}_{\min}$		& Maximal negative pressure loss allowed in the network [mbar] \\

			Cost parameters			& $\alpha_\text{p}^\text{e}$	& Before--tax unit price of electricity [\euro/kWh] \\
			\:						& $\alpha_\text{p}^\text{g}$	& Before--tax unit price of gas [\euro/kWh] \\
			\:						& $\beta^\text{e}$	& Tax rate for electricity [\euro/kWh] \\
			\:						& $\beta^\text{g}$	& Tax rate for gas [\euro/kWh] \\
			\:						& $t_\text{adv}$	& Tax advantage factor for electricity self-production \\
			\:						& $\alpha_\text{a}^\text{e}$	& Allocation unit price of electricity [\euro/kWa] \\
			\:						& $\alpha_\text{a}^\text{g}$	& Allocation unit price of gas [\euro/kWa] \\
			\:						& $\zeta_{i,j}^\text{g}$		& Price for installing a gas pipe on arc $(i,j)$ [\euro/a] \\
			\:						& $\gamma^t$		& Investment and maintenance cost for technology $t$ [\euro/a]\\
			\:						& $\nu_1$			& Unit cost of first stage renovation [\euro/kWh] \\
			\:						& $\nu_2$			& Unit cost of second stage renovation [\euro/kWh] \\	
			\:						& $E_\text{target}$	& Maximal amount of CO2 emission of the network [kg/a] \\
			\:						& $\kappa^\text{e}$	& Unit CO2 emission of electricity [kg/kWh] \\
			\:						& $\kappa^\text{g}$	& Unit CO2 emission of gas [kg/kWh] \\
			\:						& $\zeta_k^\text{e}$			& Unit price of installing cable type $k$ [\euro/m] \\

			\\
			\hline
		\end{longtable} 

\section*{Acknowledgments}
\label{section:Acknowledgments}
The first author is funded by the German Federal Environmental Foundation (Deutsche Bundesstiftung Umwelt -- DBU). The third author is partially financed by the German Research Foundation (Deutsche Forschungsgemeinschaft -- DFG) via the German Excellence Strategy.
Furthermore, we want to thank Dr.\ D.\ König and Dr.\ J.\ Lu for answering our questions as well as giving valuable comments on this work.



\printbibliography

\end{document}